\def\phi{\varphi}

\def\su{\subseteq}
\def\a{\alpha}
\def\b{\beta}
\def\g{\gamma}

\def\l{\lambda}
\def\k{\kappa}
\def\z{\zeta}

\def\om{\omega}

\def\lng{\langle}
\def\rng{\rangle}
\def\ov{\overline}

\def\cf{{\operatorname {cf}}}

\def\On{{\text {On}}}

\documentclass[12pt,reqno]{amsart}
\usepackage{amstext}
\usepackage{amssymb}
\usepackage{amsfonts}

\newtheorem{lemma}{Lemma}
\newtheorem{theorem}[lemma]{Theorem}

 \newtheorem{claim}[lemma]{Claim}
 \newtheorem{fact}[lemma]{Fact}

\begin{document}
\title[trichotomy theorem]{The Condition in the Trichotomy 
Theorem is Optimal} 
\author{Menachem Kojman}
\address{Department of Mathematics and Computer Science\\
Ben Gurion University of the Negev\\
Beer Sheva, Israel}
\address{Department of Mathematical Sciences\\
Carnegie-Mellon University\\
Pittsburgh, PA}

\email{kojman@cs.bgu.ac.il}
\thanks{The first author was partially supported by  NSF grant 
No. DMS-9622579}

\author{Saharon Shelah}
\address{Institute of Mathematics\\
The Hebrew University of Jerusalem\\
Jerusalem 91904, Israel} 

\email{shelah@math.huji.ac.il} 

\thanks{The second author was partially supported by the Israeli 
Foundation for Basic Science.  Number 673 in list of publications}

\begin{abstract}
We show that the assumption $\l>\k^+$ in the Trichotomy Theorem cannot 
be relaxed to $\l>\k$.  
% This answers a question of Renling Jin.	
\end{abstract}
\maketitle

%\maketitle
\section{Introduction}

The Trichotomy Theorem specifies three alternatives for the structure 
of an increasing sequence of ordinal functions modulo an ideal on an 
infinite cardinal $\k$ --- provided the sequence has regular length 
$\l$ and $\l$ is at least $\k^{++}$.

The natural context of the Trichotomy Theorem is, of course, pcf 
theory, where a sequence of ordinal functions on $\k$ usually has 
length which is larger than $\k^{+\k}$.  However, the trichotomy 
theorem has already been applied in several proofs to sequences of length $\k^{+n}$, 
($n\ge 2$) (see  \cite{MgSh:204}, \cite{James} and \cite{eub}).

Therefore, a natural question to ask is, whether the Trichotomy 
Theorem is valid also for sequences of length $\k^+$, namely, whether 
the lower bound on the length of the sequence can be lowered by one 
cardinal.

Below we show that  the assumption $\l\ge\k^{++}$  in 
the Trichotomy Theorem is tight.  For every infinite $\k$, we construct an 
ideal $I$ over $\k$ and $<_I$-increasing sequence $\ov f\su \On^\k$ so 
that all three alternatives in the Trichotomy theorem are violated by 
$\ov f$.

% The idea of the construction is quite simple, and we describe it now. 
% Two sequences of functions are constructed: $\lng f_\a:\a<\om_1\rng$ 
% and $\lng g_\a:\a<\om_1\rng$ so that the values of both $g_\a$ and 
% of $f_\a$ come from some set $A_{\a+1}$ is a filration of $\omega_1$. 
% An ideal is organized so that modulo the ideal $g_\a$ is the immediate 
% successor of $f_\a$ is the reduced product of $A_{\b+1}$. Furthermore, 
% modulo the ideal $\lng f_\a:\a<\omega_1\rng$ is increasing and $\lng 
% g_\a:\a<\om_1\rng$ is decreasing. Now there cannot be a exac upper 
% bound to $\ov f$, as such a bound would have to\ldots

\section{The Counter-example}

Let $\k$ be an infinite cardinal. Denote by $\On^\k$ the class of all 
functions from $\k$ to the ordinal numbers. 

Let $I$ be an ideal over $\k$.  We write $f<_I g$, for $f,g\in \On^k$, 
if $\{i<\k:f(i)\ge g(i)\}\in I$ and we write $f\le_I g$ if 
$\{i<\k:f(i)>g(i)\}\in I$.  A sequence $\ov f=\lng f_\a:\a<\l\rng\su 
\On^\k$ is \emph{$<_I$-increasing} if $\a<\b<\l$ implies that $f_\a<_I 
f_\b$ and is \emph{$<_I$-decreasing} if $\a<\b<\l$ implies that $f_\b<_I 
f_\a$.

A function $f\in \On^k$ is a least upper bound mod $I$ of a sequence $\ov 
f=\lng f_\a:\a<\l\rng\su \On^\k$ if $f_\a<_I f$ for all $\a<\l$ and 
whenever $f_\a<_I g$ for all $\a$ then $f\le_I g$.  A function $f\in 
\On ^\k$ is an \emph{exact upper bound} of $\ov f$ if $f_\a\le f$ for 
all $\a<\l$, and whenever $g<_I f$, there exists $\a<\l$ such that 
$g<_I f_\a$. For subsets $t,s$ of $\k$, write $t\su_I s$ if $s-t\in 
I$. 
 
The dual filter $I^*$ of an ideal $I$ over $\k$ is the set of all 
complements of members of $I$.  The relations $\le_I$, $<_I$ and 
$\su_I$ will also be written as $\le_{I^*}$, $<_{I^*}$ and 
$\su_{I^*}$.

 Let us quote the theorem under discussion:

\begin{theorem} \label{trichotomy}(The Trichotomy Theorem)

Suppose $\l\ge\k^{++}$ is regular, $I$ is an ideal over $k$ and $\ov f= \lng
f_\a:\a<\l\rng$ is a $<_I$-increasing sequence of ordinal functions
on $k$. Then $\ov f$ satisfies one of the following conditions:
\begin{enumerate}
\item  $\ov f$ has an exact upper bound $f$ with $\cf f(i)>\k$ for
all $i<\k$;

\item   there are sets $S(i)$ for $i<\k$ satisfying $|S(i)|\le
\k$ and an ultra-filter $U$ over $k$ extending the dual of $I$ so that
for all $\a<\l$ there exists $h_\a\in\prod_{i<\k} S(i)$ and $\b<\l$ such that
$f_\a <_U h_\a <_U f_\b$.

\item  there is a function $g:\k\to \On$ such that the sequence 
$\ov t=\lng t_\a:\a<\l\rng$ does not stabilize modulo $I$, where
$t_\a=\{i<\k: f_\a(i)>g(i)\}$.
% %  (notice that $\ov t$ is
% $\su_I$-increasing, because $\ov f$ is $<_I$-increasing).
\end{enumerate}
\end{theorem}

Proofs of the Trichotomy Theorem are found in \cite{CA}, II,1.2, in 
\cite{eub} or in the future version of \cite{ABC}.

\begin{theorem} \label{counterEx}
For every infinite $\k$ there exists an ultrafilter $U$ over $\k$ and 
a $<_U$-increasing sequence $\ov f=\lng f_\a:\a<\k^+\rng\su \On^\k$ such 
conditions 1, 2 and 3 in the Trichotomy Theorem fail for $\ov f$.
\end{theorem}

% We observe first that a sequence as in the theorem above violates all 
% 3 possibilities in the Trichotomy. The third possibility is violated 
% because $U$ is an ultrafilter. Thus, ir remains only to prove 
% \ref{cpinterEx}.

\begin{proof}
Let $\l=\k^+$.  

% We shall shortly define a sequence $\ov f=\lng 
% f_\a:\a<\l\rng$ of ordinal functions on $\k$, but we make a few 
% preparations first.
% 
% For every $i<\k$ we  define below an increasing sequence of sets 
% $A^i_\a$ of ordinals for $\a<\l$, and then require in the definition 
% of   $\ov f$ that $f_\a(i)\in A_{\a+1}^i$. The ideal over $\k$ will 
% be chosen last. 

Let us establish some notation.

We recall that every ordinal has an expansion in base $\l$, namely can 
be written as a unique finite sum $\sum_{k\le l}\l^{\b_l}\a_l$ so that 
$\b_{k+1}<\b_k$ and $\a_k<\l$.  We limit ourselves from now on to 
ordinalz $\z<\l^\om$.  For such ordinals, the expansion in base $\k$ 
contains only finite powers of $\l$ (that is, every $\b_k$ is a 
natural number).  

We agree to write an ordinal $\z=\l^l\a_l + \l^{l -1}\a_{l-1} + \dots 
+\a_0$ simply as a finite sequence $\a_l\a_{l-1}\dots\a_0$.  We 
identify an expansion with $\l$ digits with one with $n>l$ digits by 
adding zeroes on the left.  If $\z=\a_l\a_{l-1}\ldots\a_0$, we call 
$\a_k$, for $k\le l$, the \emph{$k$-th digit} in the expansion of 
$\z$.

For $\a<\l$ and an integer $l$, define:

\begin{gather}
	A^l_\a=\{\a_l\a_{l-1}\ldots\a_0:\a_k<\a \mbox{ for all } k\le l\}
\end{gather}
 
$A_\a^l$ is the set of all ordinals below ${\l}^\omega$ whose 
expansion in base $\l$ contains $l+1$ or fewer digits from $\a$.  

\begin{fact}\label{fact}
For all $\a<\l$ and $l<\om$,
		\begin{enumerate}  
		\item $\bigcup_{\a<\l}A^l_\a=\l^{l+1}$
		
		\item The ordinal $\sum_{k=0}^{l}\l^k=
		\overbrace{\a\a\ldots\a}^{l+1}$ is the maximal element in 
		$A^l_{\a+1}$
		
		 \item if $\z=\a_l\a_{l-1}\ldots\a_0\in A^l_\z$ is not maximal 
		 in $A^l_{\a+1}$, then the immediate successor of $\z$ in $A^n_{\a+1}$ 
		 is obtained from $\z$ as follows: let $k$ be the first $k\le 
		 l$ for which $\a_k<\a$.  Replace $\a_k$ by $\a_k+1$  and 
		 replace $\a_m$ by $0$ for all $m<k$
	\end{enumerate}
\end{fact}

Fix a partition $\{X_n:n<\om\}$  of $\k$ with $|X_n|=\k$ for all 
$n$. Let $n(i)$, for $i<\k$, be the unique $n$ for which $i\in X_n$.

By induction on $\a<\k^+$, define $f_\a:\k\to \On$ so that:

\begin{enumerate}
	\item  $f_\a(i)\in A_{\a+1}^{n(i)}-A_\a^{n(i)}$ 

	\item For all $n,l < \om$ and  finite, strictly increasing, 
	sequences $\lng \a_k:k\le l\rng\su\l$ it holds that for every sequence $\lng 
	\z_k:k\le l\rng$ which satisfies  $\z_k\in 
	A^n_{\a_k+1}-A_{\a_k}$,  there are $\k$ many $i\in X_n$ for 
	which $\bigwedge_{k\le l}f_{\a_k}(i)=\z_k$.
 \end{enumerate}

The first item above says that $f_\a(i)$ is an ordinal below $\l^\om$ 
whose expansion in base $\l$ has $\le n(i)$ digits, at least one of 
which is $\a$.  The second item says that every possible finite 
sequence of values $\lng \z_k:k<l\rng$ is realized $\k$ many times as 
$\lng f_{\a_\k}(i):k<\l\rng$ for an arbitrary increasing sequence 
$\lng \a_\k:k<l\rng$.

  The induction required to define the sequence is straightforward.

Define now, for every $\a<\k^+$, a function 
$g_\a:\k\to \On$  as follows:

\begin{gather}\label{g}
	g_\a(i)=\min[(A_{\a+1}^n(i)\cup\{\l+^{n(i)+1}\}) - f_\a(i)]
\end{gather}

Since $f_\a(i)<\l^{n(i)+1}$ for $i\in X_n$, the definition is good.  
If $f_\a(i)$ is not maximal in $A^{n(i)}_{\a+1}$,  then $g_\a(i)$ is 
the \emph{immediate successor of $f_\a(i)$ in $A^{n(i)}_{\a+1}$}. Let 
us make a note of that:

\begin{fact}\label{note}
	There are no members of $A_{\a+1}^{n(i)}$ between $f_\a(i)$ and 
	$g_\a(i)$
\end{fact}

We have  defined two sequences:

\begin{gather}
	\ov f=\lng f_\a:\a<\l\rng\notag \\
    \ov g=\lng g_\a:\a<\l\rng\notag 
\end{gather}
        
Next we wish to find an ideal modulo which $\ov f$ is $<_I$-increasing 
and $\ov g$ is a $<_I$-decreasing sequence of upper bounds of $\ov f$.

\begin{claim}\label{dec}
For every finite increasing sequence $\a_0<\a_1<\dots <\a_l<\l$ there 
exists $i<\k$ such that for all $k<l$

\begin{gather}\label{pat} 
f_{\a_k}(i)<f_{\a_{k+1}}(i)<g_{\a_{k+1}}(i)<g_{\a_k}(i)
\end{gather}
\end{claim}

\begin{proof} Suppose $\a_0<\a_1<\dots <\a_l<\l$ is given and choose 
$n>l$.  Let $\z_0=\overbrace{\a_0\a_0\ldots\a_0}^{l+1}\in A^n_{\a_0}$.  
Let $\z_{k+1}$ be obtained from $\z_k$ by replacing the first $l+1-k$ 
digits of $\z_\k$ by $\a_{k+1}$:

\begin{align}
	\a_0\a_1\ldots\a_k\a_{k+1}\ldots\,\a_l&=\z_l\notag\\
	\vdots\hskip35pt&{}\notag\\
	\a_0\a_1\ldots\a_{k-1}\a_k\ldots\a_k&=\z_k\notag\\
	\vdots \hskip35pt&{}\notag\\
	\a_0\a_1\;\ldots\ldots\;\a_1\ldots\,\a_1&=\z_1\notag\\
	\a_0\a_0\;\ldots\ldots\;\a_0\ldots\,\a_0&=\z_0\notag
\end{align}

Thus $\z_0<\z_1<\ldots<\z_l$ and $\z_k\in A^{l}_{\a_k}\su A^n_{\a_k}$ is 
not maximal in $A^n_{\a_k}$ (because $n>l$).  Let $\xi_k$ be the 
immediate successor of $\z_k$ in $A^n_{\b_k}$.

By  Fact \ref{fact} above,  we have

\begin{align}
	1\overbrace{0\;\ldots\quad\,\ldots\;\,0\,\ldots\quad\ldots\;0}^{l+1}&=\xi_0\notag\\
	(\a_0+1)\,0\;\ldots \ldots\,\ldots\ldots0&=\xi_1\notag\\	
	\vdots\hskip35pt&{}\notag\\
	\a_0\a_1\ldots(\a_{k-1}+1)0\ldots0&=\xi_k\notag\\
    \vdots \hskip35pt&{}\notag\\
	\a_0\a_1\;\ldots\,\ldots\;(\a_{l-1}+1)0&=\xi_l\notag
\end{align}

Therefore $\z_0<\z_1<\ldots\z_l<\xi_l<\xi_{l-1}<\ldots<\xi_0$.  To 
complete the proof it remains to find some $i\in X_n$ for which 
$f_{\a_k}(i)=\z_\k$ for $k\le l$, and, consequently, by the definition 
(\ref{g}) above, $g_{\a_k}(i)=\xi_k$.  The existence of such $i\in 
X_n$ is guaranteed by the second condition in the definition of $\ov 
f$.
\end{proof}

 For $\a<\b<\l$, let
 
 \begin{gather}
 	C_{\a,\b}=\{i<\k: f_\a(i)<f_\b(i)<g_\b(i)<g_\a(i)\}
 \end{gather}

   \begin{claim}
 	$\{C_{\a,\b}:\a<\b<\k^+\}$ has the finite intersection property
    \end{claim}
    
     \begin{proof} 
Suppose that $\a_0,\b_0,\a_1,\b_1,\dots,\a_l,\b_l$ are given and 
$\a_k<\b_k<\l$ for $k\le l$.  Let $\lng \gamma_m:m<m(*)\rng$ be the 
increasing enumeration of $\bigcup_{k\le l}\{\a_k,\b_k\}$. To show 
that $\bigcap_{k\le l}C_{\a_k,\b_k}$ is not empty, it  
suffices to find some $i<\k$ for which the sequence $g_{\g_m}(i)$ is 
decreasing in $m$ and $f_{\g_m}(i)$ is increasing in $m$.  The 
existence of such an $i<\k$ follows from Claim \ref{dec}.
  \end{proof}
  
Let $U$ be any ultrafilter extending $\{C_{\a,\b}:\a<\b<\l\}$.  Since 
for every $\a<\b$ it holds that $f_\a<_U f_\b <_U g_\b <_U g_\a$, we 
conclude that $\ov f$ is $<_U$-increasing, that $\ov g$ is 
$<_U$-decreasing and that every $g_\a$ is an upper bound of $\ov f$ 
mod $U$.

\begin{claim}
	There is no exact upper bound of $\ov f$ mod $U$.
\end{claim}

\begin{proof}
It suffices to check that there is no $h\in \On^\k$ that satisfies 
$f_\a<_U h<_U g_\a$ for all $\a<\k^+$.  Suppose, then, that $h\in\On^k$ 
satisfies this.  Since $h<_U g_0$, we may assume that  
$g(i)<\l^{n(i)+1}$ for all $i<\k$ (by changing $h$ on a set 
outside of $U$).  

Let $i<\k$ be arbitrary.  Since $\bigcup_{\a<\l} 
A^{n(i)}_\a=\l^{n(i)+1}$, there is some $\a(i)$ so that $h(i)\in 
A^n_\a(i)$.  By regularity of $\l$ it follows that there is some 
$\a(*)<\l$ such that $h(i)\in A^{n(i)}_{\a(*)}$ for all $i<\k$.  By 
our assumption about $h$, $f_{\a(*)}<_U h <_U g_{\a(*)}$.  Thus, there 
is some $i<\k$ for which $f_{\b(*)}(i)<h(i)<g_{\b(*)}(i)$.  However, 
all three values belong to $A^{n(i)}_{\a(*)+1}$, while by Fact 
\ref{note} there are no members of $A^{n(i)}_{\a(*)+1}$ between 
$f_{\beta(*)(i)}$ and $g_{\a(*)}(i)$ --- a contradiction.
\end{proof}

 \begin{claim}
 	there are no sets $S(i)\su \On$ for $i<\k$ which satisfy condition 2 
 	in the trichotomy for $\ov f$ and $U$.
 \end{claim}

\begin{proof}
	Suppose that $S(i)$, for $i<\k$, and $h_\a\in \prod_{i<\k}S(i)$ 
	satisfy 2.  in the Trichotomy Theorem.  Find $\a<\l$ such that 
	$S(i)\su A^{n(i)}_\a$ for all $i$.  Thus $f_\a<_U h_\a<_U g_\a$ 
	--- contradiction to \ref{note}.
\end{proof} 

\begin{claim}
	there is no $g:\k\to \On$ such that $g,\ov f$ and the dual of $U$ satisfy 
	condition 3. in the Trichotomy Theorem.
\end{claim}

\begin{proof}
	Let $g:\k\to On$ be arbitrary, and let 
	$t_\a=\{i<\k:f_\a(i)>g(i)\}$.  As $\ov f$ is $<_U$-increasing, for 
	every $\a<\b<\l$ necessarily $t_\a\su_U t_\b$.  Since $U$ is an 
	ultrafilter, every $\su_U$-increasing sequence of sets stabilizes.
\end{proof}
\end{proof}


\begin{thebibliography}{1}

\bibitem{James}
James Cummings.
\newblock Collapsing successors of singulars.
\newblock {\em Proc. Amer. Math. Soc.}, 125(9):2703--2709, 1997.

\bibitem{ABC}
Menachem Kojman.
\newblock The abc of pcf.
\newblock {\em Logic Eprints}.

\bibitem{eub}
Menachem Kojman.
\newblock Exact upper bounds and their uses in set theory.
\newblock {\em submitted}, 1997.

\bibitem{MgSh:204}
Menachem Magidor and Saharon Shelah.
\newblock {When does almost free imply free? (For groups, transversal etc.)}.
\newblock {\em {Jour. Amer. Math. Soc.}},  7:769--830, 1994.

\bibitem{CA}
Saharon Shelah.
\newblock {\em Cardinal Arithmetics}, volume~29 of {\em Oxford Logic Guides}.
\newblock Oxford Science Publications, 1994.

\end{thebibliography}
\end{document}